\documentclass[11pt]{article}
\usepackage{amssymb,epsfig}

\def\<{\langle}
\def\>{\rangle}
\def\pf{\medskip\noindent{\it Proof.}~}
\def\C{\mathbb{C}}
\def\Z{\mathbb{Z}}

\def\a{\alpha}
\def\b{\beta}
\def\g{\gamma}
\def\l{\lambda}
\def\G{\Gamma}
\def\sm{\smallsetminus}

\newtheorem{thm}{Theorem}[section]
\newtheorem{lem}[thm]{Lemma}
\newtheorem{cor}[thm]{Corollary}
\newtheorem{conj}{Conjecture}

\title{Generalised triangle groups of type $(3,3,2)$}
\author{James Howie\\
Department of Mathematics and\\ Maxwell Institute for Mathematical Sciences\\
Heriot-Watt University\\
Edinburgh EH14 4AS\\
UK}

\begin{document}

\maketitle

\begin{abstract}
If $G$ is a group with a presentation of the form
$\< x,y|x^3=y^3=W(x,y)^2=1\>$, then either $G$ is
virtually soluble or $G$ contains a free subgroup of rank $2$.
This provides additional evidence in favour of a conjecture
of Rosenberger.
\end{abstract}

\section{Introduction}

A {\em generalised triangle group} is a group $G$ with a presentation of
the form $$\< x,y|x^p=y^q=W(x,y)^r=1\>$$ where $p,q,r\ge 2$ are integers
and $W(x,y)$ is a word of the form $x^{\a(1)}y^{\b(1)}\cdots x^{\a(k)}y^{\b(k)}$,
($0<\a(i)<p$, $0<\b(i)<q$).  We say that $G$ is of {\em type} $(p,q,r)$.
The parameter $k$ is called the {\em length}.  Without loss of generality,
we assume that $p\le q$.

A conjecture of Rosenberger \cite{Ros} asserts that a Tits alternative holds for
generalised triangle groups:

\begin{conj}[Rosenberger]\label{RC}
Let $G$ be a generalised triangle group.  Then
either $G$ is soluble-by-finite
or $G$ contains a non-abelian free subgroup.
\end{conj}

This conjecture
has been verified in a large number of special cases.
(See for example the survey in \cite{FRR}.) In particular it is now known:
\begin{itemize}
\item when $r\ge 3$ \cite{FLR};
\item when $\frac1{p}+\frac1{q}\ge\frac12$ \cite{BMS,Howie};
\item when $q\ge 6$ \cite{LR,Will2,BK,BK2,BKB,BKB2,HW};
\item when $k\le 6$ \cite{Ros,LR,Will};
\item for $(p,q,r)=(3,4,2)$ \cite{BKB1,HW2};
\item for $(p,q,r)=(2,4,2)$ and $k$ odd \cite{BK3}.
\end{itemize}

In the present note we describe a proof of the Rosenberger Conjecture
for the case $(p,q,r)=(3,3,2)$.  Using essentially the same argument,
we also prove the Conjecture in the case where $(p,q,r)=(2,3,2)$ and $k$
is even -- with the exception of $6$ groups that our methods are unable
to handle.

The proofs rely to some extent on computations
using the computer algebra package GAP \cite{GAP}.
The strategy of proof, however, is straightforward.  Firstly, a
theoretical analysis shows that, if $G$ is a generalised triangle
group $G$ of type $(3,3,2)$ that does not contain a non-abelian free subgroup,
then the corresponding {\em trace polynomial} must have a very restricted
form.  In particular this analysis provides a bound $k\le 20$
for the length parameter $k$ of such a group.
 Secondly, a computer search finds all words of
length up to $20$  whose trace polynomial has this restricted form.
There turn out to be precisely $19$ such words (up to a standard
equivalence relation), of which $8$
have length $k\le 6$: the conjecture is already known for the
$8$ groups corresponding to these short words.
Finally, it is observed that, in the remaining $11$ cases, a small
cancellation condition applies to $G$ (regarded as a quotient
of $\Z_3\ast\Z_3$).  We complete the proof by showing that the small
cancellation condition implies the existence of a non-abelian
free subgroup.  The small cancellation arguments applied to do this
yield somewhat more general results which may be of independent
interest, so we present these in a more general form in \S\,\ref{sc} below.

\medskip
The theoretical analysis of the case $(p,q,r)=(2,3,2)$, $k$ even,
is identical, subject to two provisos.  Firstly, an equivalence class
of words in the
$(3,3,2)$ case can correspond to either one or two equivalence classes
of words of even length in the $(2,3,2)$ case. (Here {\em equivalence}
refers to some standard moves on words $W$ that do not change the
trace polynomial or the isomorphism type of the resulting group.
See \S \ref{equiv} for details.)
Secondly, these
words are twice as long as their $(3,3,2)$ counterparts.  The latter
difference means that fewer of them are already dealt with by existing
results.

Section \ref{sc} below contains the small-cancellation results mentioned
above.  Sections \ref{equiv} and \ref{tp} contain respectively a discussion of the
equivalence relation on words, and some elementary results on their trace polynomials.
The proof of the main result on generalised triangle groups of type $(3,3,2)$ is in
Section \ref{s332}, and in Section \ref{s232} we discuss the variations needed to
attack generalised triangle groups of type $(2,3,2)$ with even length
parameter. Section \ref{comp} contains some remarks about the computational
aspects of the work, including a description of the search algorithm.
 Logs of GAP sessions performing some of the calculations are contained
in an Appendix. Tables at the end of the paper list all words
(up to equivalence) whose trace polynomials do not immediately imply the
existence of free subgroups in the corresponding group.
An ancillary file attached to this preprint contains the GAP code listings
used in the search algorithm.

\bigskip\noindent{\bf Acknowledgement}

I am grateful to Gerald Williams
for helpful discussions on some of the work presented in this paper.

\section{Small Cancellation}\label{sc}

In this section we prove two results on one-relator products of groups
where the relator satisfies a certain small cancellation condition.
We will apply these specifically to generalised triangle groups of
types $(3,3,2)$ and $(2,3,2)$ respectively, but as the results seem
of independent interest, we prove them in the widest generality available.

Suppose that $\G_1,\G_2$ are groups, and $U\in\G_1\ast\G_2$ is a cyclically reduced
word of length at least $2$. (Here and throughout this section, {\em length}
means length in the free product sense.)
A word $V\in\G_1\ast\G_2$ is called a {\em piece}
if there are words $V',V''$ with $V'\ne V''$, such that each of 
$V\cdot V'$, $V\cdot V''$ is cyclically reduced as written, and each is
equal to a cyclic conjugate of $U$ or of $U^{-1}$.  A cyclic subword
of $U$ is a {\em non-piece} if it is not a piece.

By a {\em one-relator product} $(\G_1\ast\G_2)/U$ of groups $\G_1,\G_2$
we mean the quotient of their free product $\G_1\ast\G_2$ by the normal
closure of a cyclically reduced word $U$ of positive length.
Recall \cite{H1} that a {\em picture} over the one-relator product
$G=(\G_1\ast\G_2)/U$ is a graph $\mathcal{P}$ on a surface $\Sigma$ (which
for our purposes will always be a disc) whose corners are labelled
by elements of $\G_1\cup\G_2$, such that
\begin{enumerate}
 \item the label around any vertex, read in clockwise order, spells out a cyclic
permutation of $U$ or $U^{-1}$;
\item the labels in any region of $\Sigma\smallsetminus\mathcal{P}$
either all belong to $\G_1$ or all belong to $\G_2$;
\item if a region has $k$ boundary components labelled by words $W_1,\dots,W_k\in\G_i$
(read in anti-clockwise order; with $i=1,2$), then the quadratic equation
$$\prod_{j=1}^k X_jW_jX_j^{-1}=1$$
is solvable for $X_1,\dots,X_k$ in $\G_i$.
(In particular, if $k=1$ then $W_1=1$ in $\G_i$).
\end{enumerate}

Note that edges of $\mathcal{P}$ may join vertices to vertices,
or vertices to the boundary $\partial\Sigma$,
or $\partial\Sigma$ to itself, or may be simple closed curves disjoint from the
rest of $\mathcal{P}$ and from $\partial\Sigma$.

The {\em boundary label} of $\mathcal{P}$ is the product of the labels around $\partial\Sigma$.
By a version of van Kampen's Lemma, there is a picture with boundary label
$W\in\G_1\ast\G_2$ if and only if $W$ belongs to the normal closure of $U$.

A picture is {\em minimal} if it has the fewest possible vertices among
all pictures with the same (or conjugate) boundary labels.  In particular
every minimal picture is
{\em reduced}: no edge $e$ joins two distinct vertices in such a way that the labels
of these two vertices that start and finish at the endpoints of $e$
are mutually inverse.

In a reduced picture, any collection of parallel edges between two vertices
(or from one vertex to itself)
corresponds to a collection of consecutive $2$-gonal regions, and the labels within these
$2$-gonal regions spell out a piece:

\begin{center}
\scalebox{0.4}[0.4]{\includegraphics{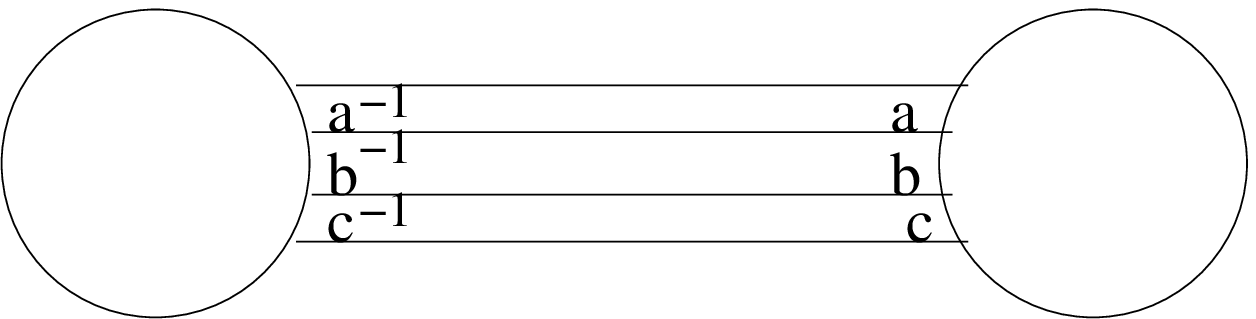}}
\end{center}

Since $U$ is cyclically reduced, no corner of an interior vertex is contained in a $1$-gonal region.

\begin{thm}
Let $\ell$ be an even positive integer.
Suppose that $U\equiv U_1\cdot U_2\cdot U_3\cdot U_4\cdot U_5\cdot U_6\in\G_1\ast\G_2$
with each $U_i$ a non-piece of length at least $\ell$.  Suppose also that $A,B\in\G_1\ast\G_2$
are reduced words of length $\ell$ such that $A$ is not equal to any cyclic conjugate
of $B^{\pm 1}$ and such that no $U_i$
is equal to a subword of a power of $A$.  Then $G:=(\G_1\ast\G_2)/\<\<U\>\>$ contains a non-abelian
free subgroup.
\end{thm}

\pf
Since $\ell$ is even and positive, any reduced word of length $\ell$ in
$\G_1\ast\G_2$ is cyclically
reduced.  Replacing $A$ by $A^{-1}$ and/or $B$ by $B^{-1}$ if necessary,
we may assume that each of $A,B$ begins with a letter from $\G_1$ and ends
with a letter from $\G_2$.
Choose a large positive integer $N>20K\ell$, where $K$ is the length of $U$,
and define $X:=A^NB^N$, $Y:=B^NA^N$.  We claim that $X,Y$ freely generate
a free subgroup of $G$.

We prove this claim by contradiction.  Suppose that $Z(X,Y)$ is a
non-trivial reduced
word in $X,Y$ such that $Z(X,Y)=1$ in $G$.  Then there exists a
picture $\mathcal{P}$ on the disc $D^2$ over the one-relator product $G$ with boundary label $Z(X,Y)$.
Without loss of generality, we may assume that
$\mathcal{P}$ is minimal, hence reduced.

Suppose that $v$ is an interior vertex of $\mathcal{P}$.  The vertex label of $v$ is $U$ or $U^{-1}$ --
by symmetry we can assume it is $U$.  The subword $U_1$ of $U$
corresponds to a sequence of consecutive corners of $v$; at least one of these corners
does not belong
to a $2$-gonal region of $\mathcal{P}$, since $U_1$ is a non-piece.
It follows that at least one of the corners of $v$ within the subword $U_1$ of the
vertex label does not belong to a $2$-gonal region.  The same follows for
the subwords $U_2,\dots,U_6$, so $v$ has at least $6$ non-$2$-gonal corners.

\begin{center}
\scalebox{0.4}[0.4]{\includegraphics{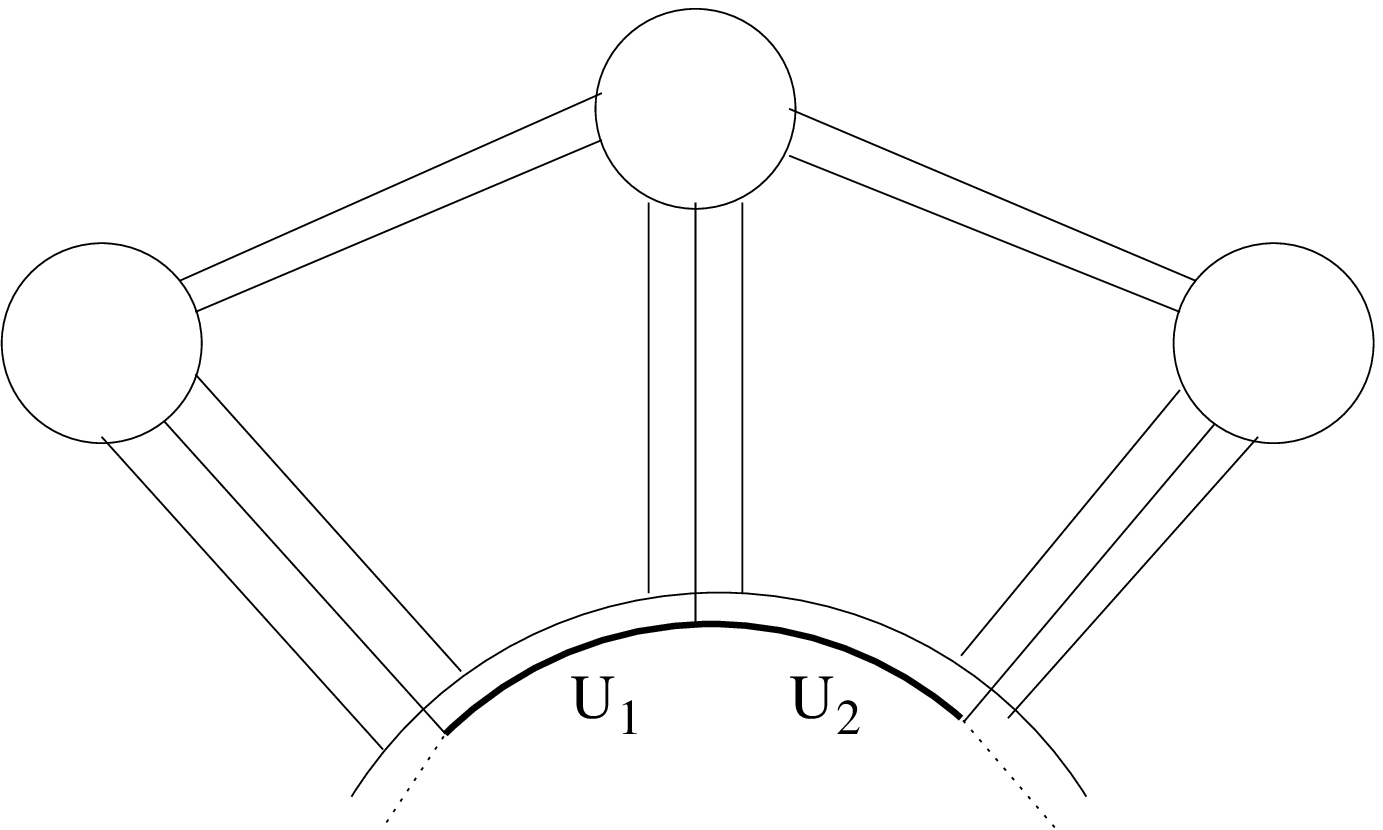}}
\end{center}

Now consider the (cyclic) sequence of {\em boundary} (that is, non-interior) vertices of
$\mathcal{P}$, $v_1,\dots,v_n$ say.  This is intended to mean that the closed path
$\partial D^2$, with an appropriate choice of starting point, meets a sequence
of arcs that go to $v_1$, separated by $2$-gons, then a sequence of arcs that go to $v_2$,
separated by $2$-gons, and so on,
finishing with a sequence of arcs that go to $v_n$, separated by $2$-gons, before returning to its starting point.
Note that it is possible that an arc of $\mathcal{P}$ joins two points on $\partial D^2$;
any such arc is ignored here.  Note also that we do not insist that $v_i\ne v_j$
for $i\ne j$ in general. It is possible for the sequence of boundary vertices to visit a vertex $v$
several times.  Nevertheless it is important to regard such visits as pairwise distinct, so the
notation $v_1,v_2,\dots$ is convenient.  We say that a boundary vertex is {\em simple} if
it appears only once in this sequence.

If $v_j$ is connected to $\partial D^2$ by $k$ arcs separated by $k-1$ $2$-gons, then
this corresponds to a common (cyclic) subword $W_j$ of $Z(X,Y)$ and $U$,
of length $k-1$.  Let $\kappa(j)\le 6$ be the maximum integer $t$ such that, for some $s\in\{1,\dots,6\}$,
$W_j$ contains a subword equal to $(U_s\cdot U_{s+1}\cdots U_{s+t})^{\pm 1}$ (indices modulo $6$).
If no such $t$ exsits, we define $\kappa(j)=-1$.

If $v_j$ is a simple boundary vertex with only $r\le 4$ corners not belonging to $2$-gons, then
it is easy to see that $\kappa(j)\ge 5-r$:

\begin{center}
\scalebox{0.4}[0.4]{\includegraphics{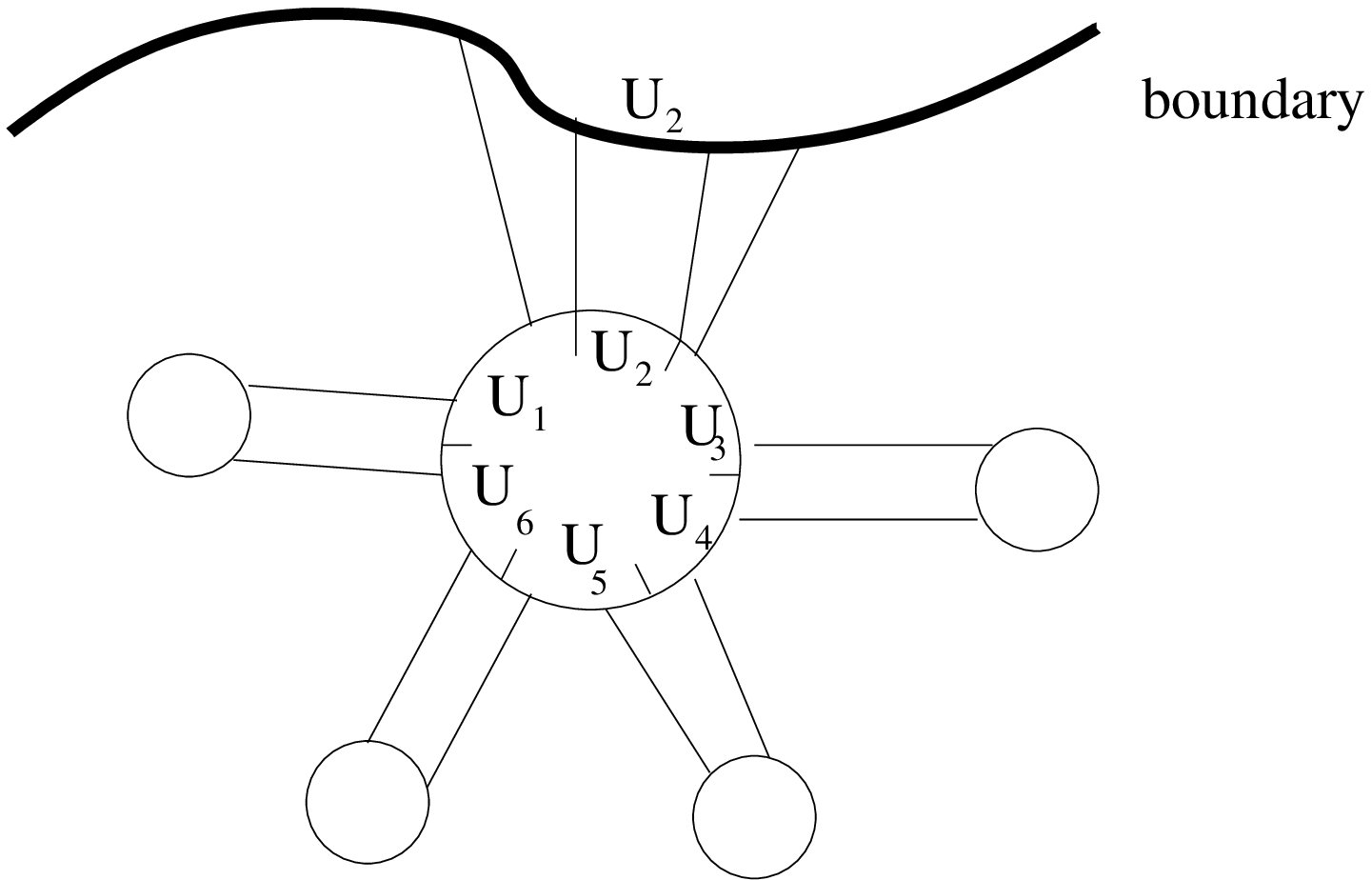}}
\end{center}

\begin{center}
\scalebox{0.4}[0.4]{\includegraphics{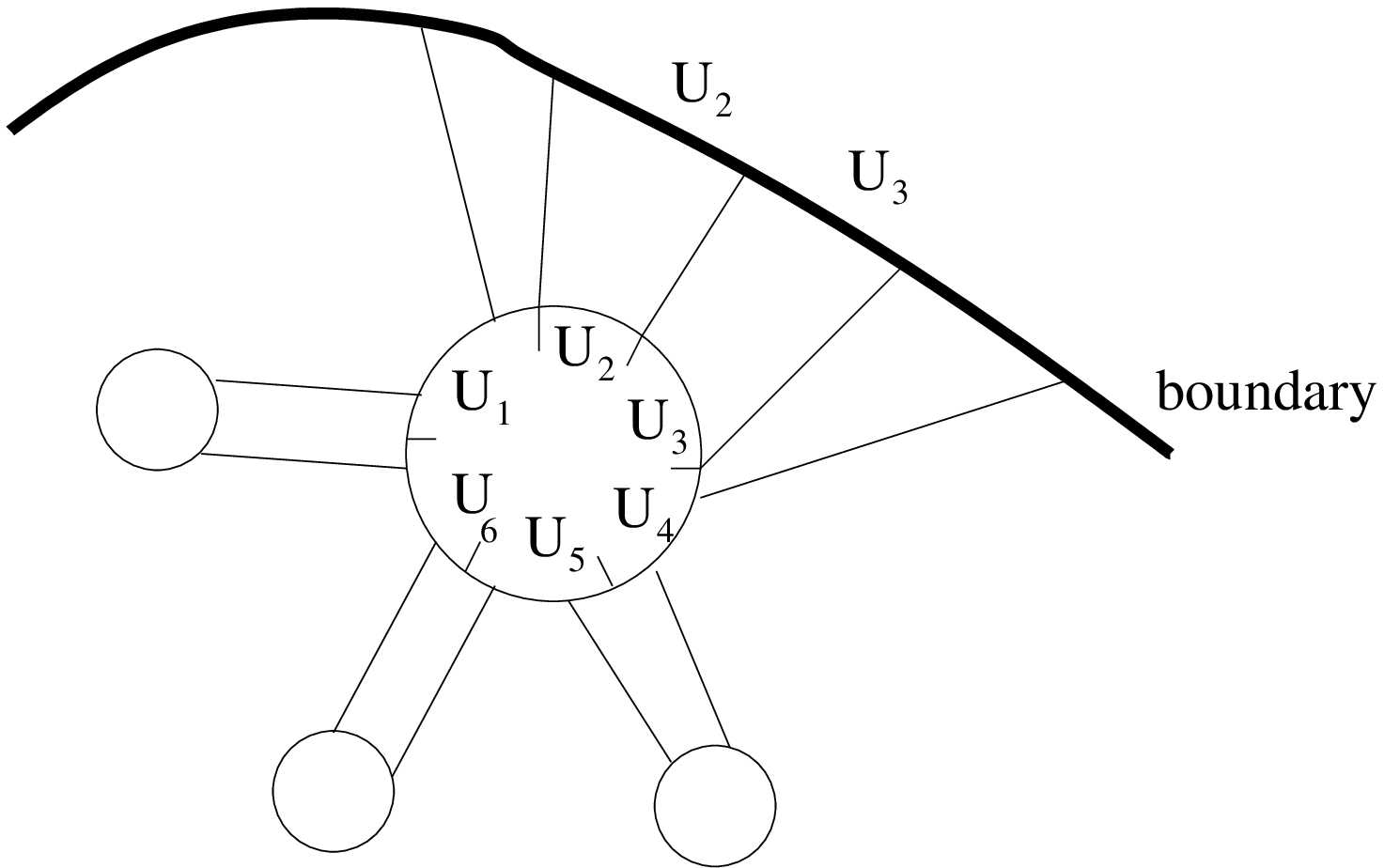}}
\end{center}

There are more complex rules for non-simple boundary vertices.  Nevertheless, it is an easy
consequence of Euler's formula, together with the fact that interior vertices have $6$ or
more non-$2$-gonal corners, that
$$\sum_{j=1}^n \kappa(j)\ge 6.$$

Now consider the word $Z(X,Y)$ as a cyclic word in $\G_1\ast\G_2$. Where a letter
$X=A^NB^N$ or $Y=B^NA^N$ is followed by another letter $X$ or $Y$, then
there is no cancellation in $\G_1\ast\G_2$.  Similarly there is no cancellation
where $X^{-1}$ or $Y^{-1}$ is followed by $X^{-1}$ or $Y^{-1}$.  Where $X$
is followed by $Y^{-1}$ or {\em vice versa}, or where $Y$ is followed by $X^{-1}$
or {\em vice versa}, then there is possible cancellation, but since $A\ne B$ the
amount of cancellation is limited to at most $\ell$ letters from either side.

If $Z$ has length $L$ as a word in $\{X^{\pm 1},Y^{\pm 1}\}$, then after cyclic reduction in $\G_1\ast\G_2$
it consists of $L$ subwords of the form $A^{\pm (N-1)}$, $L$ subwords of the form
$B^{\pm (N-1)}$, and $L$ subwords $V_1,\dots,V_L$, each of length at most $2\ell$.

Now suppose that $v_j$ is a boundary vertex of $\mathcal{P}$
with $\kappa(j)\ge 0$.  Then $U_i^{\pm 1}$ is equal to a subword of $W_j$ for some $i$.
Since $U_i$ cannot be a subword of a power of $A$,
$W_j$ is not entirely contained within one of the
segments labelled $A^{\pm (N-1)}$.

If, in addition, $\kappa(j)>0$, then $W_j$ has
a subword of of the form $(U_iU_{i+1})^{\pm 1}$ (subscripts modulo $6$)
As above, $W_j$ cannot be
contained in one of the subwords $A^{\pm (N-1)}$.  If it is contained in a subword
of $B^{\pm (N-1)}$, then it is a periodic word of period $\ell$ (that is,
its $i$-th letter is equal to its ($i+\ell$)-th letter for all $i$ for which this
makes sense).  Since $U_{i+1}$ has length at least $\ell$, there are at least two
distinct subwords of $U_iU_{i+1}$ equal to $U_i$, contradicting the fact that $U_i$
is a non-piece in $U$.

Thus we see that the subwords $W_j$ of $Z(X,Y)$
corresponding to boundary vertices $v_j$ with $\kappa(j)>0$
can occur only at certain points of $Z(X,Y)$:
where an $A^{\pm (N-1)}$-segment meets a $B^{\pm (N-1)}$-segment;
or at part of one of the words $V_i$.

In particular, the number of
boundary vertices $v_j$ with $\kappa(j)>0$ is bounded above by $L(2\ell+1)$.
It follows that $$\kappa:=\sum_j \kappa(j)\le 5L(2\ell+1),$$
where the sum is taken over those boundary vertices $v_j$ with $\kappa(j)\ge 0$.

The goal is to show that the total positive contribution to the sum $\kappa$ from
those $v_j$ with $\kappa(j)>0$
is cancelled out by negative contributions to $\kappa$ from other boundary vertices.
This will show that $\kappa\le 0$, contradicting the assertion above that $\kappa\ge 6$.

Recall that $K$ is the length of $U$.  Thus each $A^{\pm (N-1)}$-segment of
$\partial\mathcal{P}$ is joined to at least $(N-1)\ell/K$ boundary vertices, at most $2$
of which (those at the ends of the segment) can make non-negative contributions to
$\kappa$.
The remaining vertices each contribute
at most $-1$ to $\kappa$. Since $N>20K\ell$, it follows that
the negative contributions outweigh the positive contributions, as required.

This gives the desired contradiction, which proves the theorem.

\begin{cor}\label{sc33}
Let $\G_1$ and $\G_2$ be groups, and suppose $x\in\G_1$ and $y\in\G_2$
are elements of order greater than $2$.  Suppose that $W\equiv U_1\cdot U_2\cdot U_3\in\G_1\ast\G_2$
with each $U_i$ a non-piece of length at least $4$.  Then $G=(\G_1\ast\G_2)/\<\<W^2\>\>$
contains a non-abelian free subgroup.
\end{cor}

\pf Let $A_1=xyxy$, $A_2=xy^{-1}xy^{-1}$, $A_3=xyxy^{-1}$ and $A_4=xyx^{-1}y^{-1}$.
Then for $i\ne j$, $A_i$ is not equal to a cyclic conjugate of $A_j^{\pm 1}$.
Hence if (say) $U_1$ is equal to a subword of a power of $A_i$,
it cannot be equal to a subword of a power of $A_j$.  Hence there is at least
one $A\in\{A_i,~1\le i\le 4\}$ with the property that no $U_i$
is equal to a subword of a power of $A$.  Now choose $B\in\{A_i,~1\le i\le 4\}\sm A$
and apply the theorem, with $U_4=U_1$, $U_5=U_2$ and $U_6=U_3$.

\begin{cor}\label{sc23}
Let $\G_1$ and $\G_2$ be groups, and suppose $x\in\G_1$ has order $2$
and $y\in\G_2$ has order greater than $2$.  Suppose that $W\equiv U_1\cdot U_2\cdot U_3\in\G_1\ast\G_2$
with each $U_i$ a non-piece of length at least $8$.  Then $G=(\G_1\ast\G_2)/\<\<W^2\>\>$
contains a non-abelian free subgroup.
\end{cor}

\pf
Let $A_1=xyxyxyxy$, $A_2=xyxy^{-1}xyxy^{-1}$, $A_3=xyxyxyxy^{-1}$ and $A_4=xyxyxy^{-1}xy^{-1}$.
As in the previous proof, we can choose $A,B\in\{A_1,A_2,A_3,A_4\}$ such that
no $U_i$ is equal to a subword of a power of $A$, and $A$ is not equal
to a cyclic conjugate of $B^{\pm 1}$, and apply the theorem, with $U_4=U_1$, $U_5=U_2$ and $U_6=U_3$.

\section{Equivalence of words}\label{equiv}

Our object of study is a group
$$G=\< x,y|x^p=y^q=W(x,y)^r=1\>$$ where
$$W(x,y)=x^{\a(1)}y^{\b(1)}\cdots x^{\a(k)}y^{\b(k)},$$
and $0<\a(i)<p$, $0<\b(i)<q$ for each $i$.

We think of the word $W$ as a cyclically reduced word in the free product
$$\Z_p\ast\Z_q=\<x,y|x^p=y^q=1\>.$$

We regard two such words $W,W'$ as {\em equivalent} if one can be transformed
to the other by moves of the following types:
\begin{itemize}
\item cyclic permutation of $W$,
\item inversion of $W$,
\item automorphism of $\Z_p$ or of $\Z_q$, and
\item (if $p=q$) interchange of $x,y$.
\end{itemize}

It is clear that, if $W,W'$ are equivalent words, then the resulting
groups
$$G=\< x,y|x^p=y^q=W(x,y)^r=1\>$$
and
$$G'=\< x,y|x^p=y^q=W'(x,y)^r=1\>$$
are isomorphic.  Hence for the purposes of studying the Rosenberger Conjecture
(Conjecture \ref{RC}) it is enough to consider words up to equivalence.

\section{Trace Polynomials}\label{tp}

Suppose that $X,Y\in SL(2,\C)$ are matrices, and $W=W(X,Y)$
is a word in $X,Y$.  Then the trace of $W$ can be calculated
as the value of a $3$-variable polynomial, where the variables
are the traces of $X$, $Y$ and $XY$ \cite{Horowitz}.  We can use
this to find and analyse {\em essential representations}
from $G$ to $PSL(2,\C)$.  (A representation of $G$ is {\em essential}
if the images of $x,y,W(x,y)$ have orders $p,q,r$ respectively.)

We can force the images $x,y$ to have orders $p,q$ in $PSL(2,\C)$
by mapping them to matrices $X,Y\in SL(2,\C)$ of trace $2\cos(\pi/p)$
and $2\cos(\pi/q)$ respectively.  Then the trace of
$W(X,Y)\in SL(2,\C)$ is given by a one-variable polynomial $\tau_W(\l)$,
where $\l$ denotes the trace of $XY$.  Since we are in practice interested
in the case where $r=2$, we obtain an essential representation by
choosing $\l$ to be a root of $\tau_W$.

We recall here some properties of $\tau_W$.  Details can be found, for
example, in \cite{FRR}.
\begin{itemize}
\item $\tau_W$ has degree $k$;
\item when $p,q\le 3$, $\tau_W(\l)$ is monic and has integer
coefficients;
\item when $p=2$, $\tau_W$ is an odd or even polynomial,
depending on the parity of $k$.
\end{itemize}

\begin{lem}\label{equivpoly}
If $p=2$, $q=3$ and $W,W'$ are equivalent, then $\tau_W(\l)=\tau_{W'}(\l)$.
\\
If $p=q=3$ and $W,W'$ are equivalent of length $k$, then either $\tau_W(\l)=\tau_{W'}(\l)$
or\\ $\tau_W(\l)=(-1)^k\tau_{W'}(1-\l)$.
\end{lem}

\pf Since the trace of a matrix is a conjugacy invariant,
it follows that the trace polynomial is unchanged by cyclically
permuting $W$.  Moreover, if $X\in SL(2,\C)$ then the traces of $X,X^{-1}$
are equal, so the trace polynomial is unchanged by inverting $W$.

Suppose first that $p=2$ and $q=3$.  Then we cannot interchange $x$ and $y$.
Moreover, there is no nontrivial automorphism of $\Z_2$ and only one
nontrivial automorphism of $\Z_3$, which replaces
$y$ by $y^2$.  If $tr(X)=0$ and $tr(Y)=1$, then $tr(Y^{-1})=1$, and
$$tr(XY^{-1})+tr(XY)=tr(X)tr(Y)=0,$$
so $tr(XY^{-1})=-tr(XY)=-\l$, so
$$tr(W(X,Y^2))=tr(W(X,-Y^{-1}))=(-1)^k\tau_W(-\l)=\tau_W(\l).$$

In other words, this does not change $\tau_W(\l)$, as claimed.

Now suppose that $p=q=3$. 
If $tr(X)=1=tr(Y)$, then $tr(Y^{-1})=1$ also. Interchanging $x,y$ in $W$
has the effect on $\tau_W(\l)=tr(W(X,Y))$ of replacing
$\l=tr(XY)$ by $tr(YX)=\l$ -- in other words, no change.

Now in this case $$tr(XY^{-1})+tr(XY)=tr(X)tr(Y)=1.$$
Hence replacing $y$ by $y^2$ has the effect of
replacing $\tau_W(\l)=tr(W(X,Y)$
by
$$tr(W(X,Y^2))=tr(W(X,-Y^{-1}))=(-1)^k tr(W(X,Y^{-1})) = (-1)^k \tau_W(1-\l),$$
as claimed.

\begin{lem}\label{332to232}
Let $W$ be a cyclically reduced word in $\Z_3\ast\Z_3=\<x,y|x^3=y^3=1\>$,
and define $Z(u,v)=W(uvu,v)\in\Z_2\ast\Z_3=\<u,v|u^2=v^3=1\>$.  Then
$\tau_Z(\lambda)=(-1)^k\tau_W(2-\l^2)$.
\end{lem}

\pf
Let $U,V$ be matrices with $tr(U)=0$, $tr(V)=1$.
Define $X=V$ and $Y=-UVU$ so that $tr(X)=1=tr(Y)$,
and $tr(XY)=-tr((UV)^2)=2-\l^2$ where $\l=tr(UV)$.
Hence
$$\tau_Z(\l)=tr(Z(U,V))=tr(W(UVU,V))=(-1)^k tr(W(X,Y))=(-1)^k\tau_W(2-\l^2)$$
as claimed

\begin{thm}\label{tpoly}
Let $G=\<x,y|x^3=y^3=W(x,y)^2=1\>$ where
$W=x^{\a(1)}y^{\b(1)}\cdots x^{\a(k)}y^{\b(k)}$ with $\a(i),\b(i)\in\{1,2\}$
for each $i$.  If $G$ does not contain a
free subgroup of rank $2$, then $\tau_W(\lambda)$
has the form
$$\tau_W(\lambda)=\lambda^a (\lambda-1)^b (\lambda^2-\lambda-1)^c$$
with $a,b\le 1$ and $c\le 3(a+b+1)$.  In particular $k=a+b+2c\le 20$.
\end{thm}

\pf If $\lambda_0$ is a root of the trace polynomial, then
there is an essential representation $\rho:G\to PSL(2,\C)$
such that $\rho(x),\rho(y)$ are represented by matrices of trace $1$
and $\rho(xy)$ is represented by a matrix of trace
$\lambda_0$.  If the image $\rho(G)$ of $\rho$ is non-elementary,
then $\rho(G)$, and hence also $G$, contains a free subgroup of rank 2,
contrary to hypothesis.

Hence every essential representation $G\to PSL(2,\C)$ has elementary
image. But the only elementary subgroups of $PSL(2,\C)$
generated by two elements of order $3$ that contain elements of order $2$
are isomorphic to $A_4$ (corresponding to roots
$0$ or $1$ of $\tau_W$)
and $A_5$ (corresponding to roots $\frac{1\pm\sqrt{5}}{2}$).

Since $\tau_W$ has integer coefficients,
the two potential roots $\frac{1\pm\sqrt{5}}{2}$ occur with
equal multiplicities. Since $p=q=3$, $\tau_W$
is monic.  Thus $\tau_W$ has the form
$$\tau_W(\lambda)=\lambda^a (\lambda-1)^b (\lambda^2-\lambda-1)^c$$
for some non-negative integers $a,b,c$.

To obtain the desired bounds on $a,b,c$ we use the following observation.
The space $\mathcal{M}_1$ of matrices in $SU(2)\subset SL_2(\C)$ with trace $1$ is path-connected.
(Indeed, it is homeomorphic to the $2$-sphere $S^2$.)
For any $X\in\mathcal{M}_1$, we can vary $Y$ continuously in $\mathcal{M}_1$
from $X$ to $X^{-1}$, and $\lambda=tr(XY)$ will vary continuously
from $-1=tr(XX)$ to $2=tr(XX^{-1})$. By the Intermediate Value Theorem,
any $\lambda\in [-1,2]$ can be realised as $tr(XY)$ for some choice of $X,Y\in\mathcal{M}_1$.
But for $X,Y\in\mathcal{M}_1$ we have $W(X,Y)\in SU(2)$, so $\tau_W(\l)=tr(W(X,Y))\in [-2,2]$.
This shows that $|\tau_W(\l)|\le 2$ for $-1\le\l\le 2$.
Now $|\tau_W(2)|=2^a$ and $|\tau_W(-1)|=2^b$, so $a\le 1$ and $b\le 1$.
Finally, $$\left|\tau_W\left(\frac{1}{2}\right)\right|=\left(\frac54\right)^c\left(\frac12\right)^{a+b}.$$
From this we deduce that $$c\ln(5)\le (a+b+2c+1)\ln(2),$$ which implies
the desired conclusion $$c\le 3(a+b+1)$$ given that $a+b\in\{0,1,2\}$.

\bigskip
Essentially the same proof gives the following parallel version:

\begin{thm}\label{tpoly2}
Let $G=\<u,v|u^2=v^3=W(u,v)^2=1\>$ where
$W=uv^{\a(1)}\cdots uv^{\a(k)}$ with $\a(i)\in\{1,2\}$
for each $i$ and $k$ even.  If $G$ does not contain a
free subgroup of rank $2$, then $\tau_W(\lambda)$
has the form
$$\tau_W(\lambda)=(\lambda^2-1)^a (\lambda^2-2)^b (\lambda^4-3\lambda^2+1)^c$$
with $a,b\le 1$ and $c\le 3(a+b+1)$.  In particular $k=2a+2b+4c\le 40$.
\end{thm}

\section{The main result}\label{s332}

\begin{thm}\label{main}
Let $G=\<x,y|x^3=y^3=W(x,y)^2=1\>$ be a generalised triangle
group of type $(3,3,2)$.  Then the Rosenberger Conjecture
holds for $G$: either $G$ is soluble-by-finite, or $G$ contains
a non-abelian free subgroup.
\end{thm}

\pf
Write $$W=x^{\a(1)}y^{\b(1)}\cdots x^{\a(k)}y^{\b(k)}.$$
A computer search
produces a list of all words $W$, up to equivalence, for which
the trace polynomial $\tau_W$ has the form indicated in Theorem
\ref{tpoly}: see Table \ref{t332}.  If $W$ is not equivalent
to a word in the list, then $G$ has a nonabelian free subgroup
by Theorem \ref{tpoly}, so we may restrict our attention to the words
$W$ in Table \ref{t332}.

For those $W$ in Table \ref{t332} for which
$k\ge 7$ (namely, numbers 9-19) the small cancellation hypotheses of Corollary \ref{sc33}
are satisfied, and so $G$ contains a nonabelian free subgroup.

For $k\le 6$ (words 1-8) in the table,
the result is known.  Specifically, groups 1-3 are well-known
to be finite
of orders $12$, $180$ and $288$ respectively;
groups 4-6 were proved to have nonabelian free subgroups
in \cite{LR}; and finally groups 7 and 8 were shown in \cite{Will}
to be {\em large}.  (That is, each contains a subgroup of
finite index which admits an epimorphism onto a non-abelian free
group.)  Since \cite{Will} is an unpublished thesis, we will
give, for each result we cite from \cite{Will}, either a GAP
calculation reproducing Williams' argument, or an independent proof.
In particular, a GAP calculation following the proof
in \cite{Will} for Group 7 is shown in the Appendix. The largeness of
Group 8 in Table \ref{t332} follows from the largness of the corresponding
group in Table \ref{t232a}, as discussed in \S\,\ref{s232} below.
The latter group is shown to be large using a separate GAP calculation,
which is also reproduced in the Appendix.

This completes the proof.

\section{Variation: type $(2,3,2)$ with even length}\label{s232}

The group $G=\<x,y|x^3=y^3=W(x,y)^2=1\>$ has a homomorphic image
$\overline{G}$ which is an index-$2$ subgroup of
$H=\<u,v|u^2=v^3=W(uvu,v)^2=1\>$.  Clearly $H$ is a generalised
triangle group of type $(2,3,2)$ whose length parameter $k$ is
twice that of $G$.  Conversely, every generalised triangle
group of type $(2,3,2)$ with even length parameter arises
in this way.  There is thus at least a superficial parallel
between generalised triangle groups of type $(2,3,2)$ with even length
parameter and those of type $(3,3,2)$.  The two types can
be analysed in entirely analogous ways.  In particular, the
same computer search used to list the possible relators in $G$
yields also the possible relators in $H$
(see Table \ref{t232a}).  However, we must take
care over a few details.

\begin{enumerate}
 \item Interchanging $y,y^2$ in $W(x,y)$ produces a word
$W'(x,y)$ that is equivalent to $W(x,y)$.  However,
$W'(uvu,v)=W(uvu,v^2)$ is not in general equivalent to $W(uvu,v)$.
Thus each candidate for $W(x,y)$ in $G$ gives rise to either one
or two candidates for $W(uvu,v)$ in $H$ (up to equivalence).
This is reflected in the numeration of words in Table \ref{t232a}:
for example, word 1, $W(x,y)=xy$ in Table \ref{t332} gives rise to words
1a, $W(uvu,v)=uvuv$, and 1b, $W(uvu,v^2)=uvuv^2$ in Table \ref{t232a}.
(Where the two words $W(uvu,v)$ and $W(uvu,v^2)$ are equivalent, only
one is shown in Table \ref{t232a}.)
\item If $H$ contains a nonabelian free subgroup, then
so does $\overline{G}$, and hence so does $G$.  We have already used this explicitly
in the proof of Theorem \ref{main}:
taking $H$ to be Group 8 in Table \ref{t232a}, we show in the Appendix
that $H$ is large.  In this case $G$ is Group 8 of Table \ref{t332}, which we also deduced
to be large.
The converse implication
does not necessarily hold, however.  This is most graphically illustrated
by the case of Group 4 in Table \ref{t332}.  As mentioned in the proof of Theorem
\ref{main}, this was shown to contain a non-abelian free subgroup in \cite{LR}.
However, one of the two corresponding groups in Table \ref{t232a}, namely Group 4a,
is known to be finite of relatively small order \cite{HMT}.
\item If $W(x,y)$ is such that $W(uvu,v)$ satisfies the small-cancellation
hypothesis of Corollary \ref{sc23},
then $W(x,y)$ satisfies the small-cancellation
hypothesis of Corollary \ref{sc33}, but the converse does not hold
in general.
\item If $W(x,y)$ has length parameter $k\in\{4,5,6\}$, then known
results imply that the Rosenberger conjecture holds for $G$, as we saw in
the proof of Theorem \ref{main}.
But $W(uvu,v)$ has length parameter $2k\in\{8,10,12\}$ and existing results
do not necessarily apply to $H$.
\end{enumerate}

These remarks indicate that the $(2,3,2)$ situation, with even length,
is somewhat more complicated than the $(3,3,2)$ case.  We have not
been able to prove the Rosenberger conjecture in its
entirety for the $(2,3,2)$ case.  Nevertheless, we have been able to reduce the
number of potential counterexamples to $6$.

\begin{thm}\label{second}
Let $H=\<u,v|u^2=v^3=Z(u,v)^2=1\>$ be a generalised triangle
group, where $Z(u,v)=uv^{\gamma(1)}\cdots uv^{\gamma(2k)}$
and $\gamma(i)\in\{1,2\}$ for each $i$.  Then the Rosenberger
conjecture holds for $H$, except possibly when $Z$ is, up to
equivalence, one of the following:
\begin{enumerate}
\item $(uv)^3(uv^2)^2uv(uv^2)^2uvuv^2$;
\item $(uv)^4(uv^2)^3(uv)^2uv^2$;
\item $(uv)^5(uv^2)^3(uv)^2uv^2uv(uv^2)^2$;
\item $(uv)^4(uv^2)^4uv(uv^2)^3(uv)^2uv^2uv(uv^2)^2$;
\item $(uv)^4(uv^2)^4uv(uv^2)^2uv(uv^2)^3(uv)^3(uv^2)^2uvuv^2$;
\item $(uv)^4(uv^2)^2uv(uv^2)^3(uv)^2uv^2uv(uv^2)^2$.
\end{enumerate}
\end{thm}

\pf
The proof of this theorem follows the same pattern as that of Theorem \ref{main}.
The same computer search that produced Table \ref{t332}
also produces a complete list (Table \ref{t232a}) of those words
(up to equivalence) whose trace polynomials have the form indicated in
Theorem \ref{tpoly2}.  If $W$ is not equivalent to a word in Table
\ref{t232a}, then $H$ contains a nonabelian free subgroup, by Theorem \ref{tpoly2}.

Table \ref{t232a} is split into three parts.  Part 3 contains the six exceptional words
listed in the statement: we can prove nothing about the corresponding groups $H$.

Each word in part 2 of Table \ref{t232a} satisfies the small-cancellation
hypothesis of Corollary \ref{sc23}, so the corresponding group $H$ contains a nonabelian
free subgroup by Corollary \ref{sc23}.

Most of the groups in part 1 of table \ref{t232a} can be handled by
existing results.  Specifically, groups 1a, 1b, 2, 3, 4a and 6 are known to be
finite of the given orders \cite{HMT,LRS}, while group 4b was shown to contain
non-abelian free subgroups in
\cite{Will} (by observing that its unique subgroup of index $2$ is the group
4 in Table 1).

The remaining two groups can be dealt with by calculations using GAP \cite{GAP}
In group 5, the normal closure of $(uv)^{10}$ has index $7680$ and is free abelian of
rank $6$, while in group 8 the normal closure of $(uv)^{5}$ has a non-abelian free
homomorphic image of rank $3$.  (Logs of GAP sessions performing these calculations are shown in the
Appendix.)

\section{Computational Aspects}\label{comp}

The main computational aspect of this work is the search for words
with appropriate trace polynomials.  By Lemma \ref{332to232}
the search in the $(3,3,2)$ case is essentially the same as that
in the $(2,3,2)$ case with $k$ even.  In what follows we use
the latter framework.  Thus we put
$$G=\<x,y|x^2=y^3=W(x,y)^2=1\>,$$
$$W(x,y)=xy^{\a(1)}\cdots xy^{\a(k)}$$
where $k$ is even and $\a(j)\in\{1,2\}$ for each $j$.

We use the formulae in \cite[Lemma 9]{HW}
for the coefficients of $\tau_W(\lambda)$ to restrict the shape of the
words for which we are searching. 
In our context, the coefficient of $\lambda^{k-2}$ in $\tau_W(\l)$,
where $W=xy^{\a(1)}\cdots xy^{\a(k)}$, is $B_1:=b(1)+\cdots+b(k)$, where
$$b(j):=\left\{\begin{array}{ll} -1\quad &\mathrm{if}~a(j)=a(j+1)\\
\frac{-1+i\sqrt{3}}{2} &\mathrm{if}~a(j)=2\ne a(j+1)\\
\frac{-1-i\sqrt{3}}{2} &\mathrm{if}~a(j)=1\ne a(j+1).\end{array}\right.$$
Moreover, the coefficient of $\l^{k-4}$ in $\tau_W(\l)$ is
$B_2:=\sum_{\{j,j'\} } b(j)b(j')$, where the sum is over all $2$-element
subsets $\{j,j'\}$ of $\{1,\dots,k\}$ such that $j\ne j'\ne j\pm 1~\mathrm{mod}~k$. 
 
Thus if we rewrite $W$
in the form
$$W(x,y)=(xy)^{\b(1)}(xy^2)^{\g(1)}\cdots(xy)^{\b(m)}(xy^2)^{\g(m)}$$
with $\b(1)+\g(1)+\cdots+\b(m)+\g(m)=k$, then it follows from \cite[Lemma 9]{HW} that
$B_1=m-k$, so
$$\tau_W(\l)=\l^k-(k-m)\l^{k-2}+\cdots$$
In particular, if $\tau_W(\l)=(\l^2-1)^a(\l^2-2)^b(\l^4-3\l^2+1)^c$ as required
by Theorem \ref{tpoly2}, then $m=a+c$.

With the above calculation in mind, it is convenient to store the word $W$ in
the form of the list $L_W:=[\b(1),\g(1),\dots,\b(m),\g(m)]$ of positive
integers.  The GAP \cite{GAP}
command `OrderedPartitions' produces all such lists.

At this point we encounter a software problem: the output of `OrderedPartitions(40,20)',
for example, should be a list of more than $2^{32}$ lists, which would
exceed GAP's upper bound for list lengths. So the command
`OrderedPartitions(40,20)' will lead to an error.  However, we
can overcome this problem as follows.

It turns out from the formulae in \cite[Lemma 9]{HW} that, just as the second coefficient $B_1$
of $\tau_W(\l)$ determines the length of the list corresponding to $W$, the
third coefficient $B_2$ determines the number of elements in that list which
are equal to $1$. To see this, note that
$$2B_2=B_1^2-\sum_{j=1}^k b(j)^2-2\sum_{j=1}^k b(j)b(j+1).$$
From $k$ and $B_1$ we can calculate $m$, and hence determine how many of the $b(j)$ are equal to
$-1$. This in turn enables us to calculate $\sum b(j)^2$.  Hence, if we also know $B_2$,
we can calculate $\sum b(j)b(j+1)$.  On the other hand, the expression 
$$W(x,y)=(xy)^{\b(1)}(xy^2)^{\g(1)}\cdots(xy)^{\b(m)}(xy^2)^{\g(m)}$$
enables us to count the number of $j$ for which $b(j)b(j+1)\ne 1$: this happens precisely
at the beginning and the end of each syllable $(xy)^{\b(i)}$ or $(xy^2)^{\g(i)}$ for which
$\b(i)\ge 2$ (resp., $\g(i)\ge 2$). Reconciling this with our previous calculation of
$\sum b(j)b(j+1)$ tells us how many entries in $L_W$ are equal to $1$.

We omit the details, but when
$$\tau_W(\l)=(\l^2-1)^a(\l^2-2)^b(\l^4-3\l^2+1)^c$$
with $a,b\le 1$ as in Theorem \ref{tpoly2}, the calculation shows that the resulting list $L_W$ has
\begin{itemize}
 \item length $2c+2$ with $c$ entries equal to $1$, if $a=b=1$;
\item length $2c+2$ with $c+1$ entries equal to $1$, if $a=1$ and $b=0$;
\item length $2c$ with $c-1$ entries equal to $1$, if $a=0$ and $b=1$;
\item length $2c$ with $c$ entries equal to $1$, if $a=b=0$.
\end{itemize}

Now an ordered partition $\mathcal{P}$ of $n$ into $p$ positive integers, of which
precisely $q$ are equal to $1$, can be completely described
by two simpler partitions.  The first is the partition of $n-p$
into $p-q$ parts, obtained by subtracting $1$ from each element of
$\mathcal{P}$ and then removing the zeroes.  The second is a partition
of $p$ into $q+1$ positive integers, which encodes which entries of
$\mathcal{P}$ are equal to $1$. This gives us
an algorithm implementable in GAP for conducting the search.
Use `OrderedPartitions' to create two lists of lists.  For each pair
in list1 $\times$ list 2, create a list that corresponds to a word.
Calculate its trace polynomial: if this matches the form in Corollary
\ref{tpoly2} then add it to the output.

In practice we refine this algorithm in a number of ways.
\begin{enumerate}
 \item We wish our output to contain only one word from each equivalence
class.  Ideally, we could select only one word from each equivalence class
before calculating the trace polynomial, but this is not very efficient
to do.  A good compromise is to apply a fast-but-crude filter before the event (which
will occassionally let through two or more equivalent words while ensuring that
at least one word from each class gets through), and then to apply a less efficient but
rigorous filter to the (much smaller) output data.
\item Before calculating trace polynomials, we apply a further filter
to check that $G$ admits the appropriate essential permutation representations
(onto $A_4$, $S_4$, $A_5$).  This reduces the number of matrix calculations
that are required.
\item We replace the trace polynomial calculation by an amended version, that
enables us to find the value of $\tau_W(\l)$ for integer values of $\l$ using only
matrices with integer entries (which is more efficient than doing matrix calculations
over a polynomial ring).  To test for correctness of the trace polynomial,
it suffices to test its values at sufficiently many integer points.
\end{enumerate}

Code implementing this algorithm is listed in the ancillary file {\em gtg232.g}.

\section*{Appendix: GAP Sessions}

1)
The following GAP session considers Group 5 from Table \ref{t232a}
$$G_5:=\< x,y|x^2=y^3=((xy)^4(xy^2)^2xyxy^2)^2=1\>.$$
It demonstrates that $G_5$ has a free abelian normal subgroup
$N$ of rank $6$ and index $7680=60\cdot 2^7$.

A more detailed analysis shows that
the kernel $K$ of the essential representation
$G_5\to A_5$ has commutator subgroup $[K,K]$ of
order $2$, and $K/[K,K]$ is free abelian of rank $6$.
It follows that the subgroup $K^2$ generated by
$\{g^2;g\in K\}$ is central in $K$ and has index $2^7$.
The subgroup $N$ above is precisely $K^2$.

\begin{verbatim}
gap> F:=FreeGroup(2);;
gap> x:=F.1;; y:=F.2;;
gap> W:=(x*y)^4 * (x*y^2)^2 * (x*y) * (x*y^2);;
gap> G:=F/[x^2,y^3,W^2];;
gap> Q:=F/[x^2,y^3,W^2,(x*y)^10];;
gap> Size(Q);
7680
gap> H:=Subgroup(G,[(G.1*G.2)^10]);;
gap> P:=PresentationNormalClosure(G,H);;
gap> SimplifyPresentation(P);
#I  there are 22 generators and 295 relators of total length 1240
#I  there are 9 generators and 77 relators of total length 380
#I  there are 6 generators and 27 relators of total length 160
#I  there are 6 generators and 22 relators of total length 114
gap> TzGoGo(P);
#I  there are 6 generators and 18 relators of total length 86
#I  there are 6 generators and 16 relators of total length 68
gap> TzPrint(P);
#I  generators: [ _x5, _x6, _x219, _x220, _x498, _x500 ]
#I  relators:
#I  1.  4  [ 3, 2, -3, -2 ]
#I  2.  4  [ -2, 5, 2, -5 ]
#I  3.  4  [ -6, -2, 6, 2 ]
#I  4.  4  [ -6, -4, 6, 4 ]
#I  5.  4  [ -4, -2, 4, 2 ]
#I  6.  4  [ -5, 3, 5, -3 ]
#I  7.  4  [ 3, 1, -3, -1 ]
#I  8.  4  [ 3, -6, -3, 6 ]
#I  9.  4  [ -3, 4, 3, -4 ]
#I  10.  4  [ 1, -6, -1, 6 ]
#I  11.  4  [ 4, 5, -4, -5 ]
#I  12.  4  [ 1, 5, -1, -5 ]
#I  13.  4  [ -5, -6, 5, 6 ]
#I  14.  4  [ 1, -2, -1, 2 ]
#I  15.  4  [ 4, -1, -4, 1 ]
#I  16.  8  [ -5, -2, 4, -3, 5, 2, -4, 3 ]
\end{verbatim}

2)
The following GAP session considers Group 8 from Table \ref{t232a}
$$G_8:=\< x,y|x^2=y^3=((xy)^4(xy^2)^3xyxy^2(xy)^2xy^2)^2=1\>.$$
It demonstrates that the kernel of the essential representation
$G_8\to A_5$ admits an epimorphism onto the free group of
rank $3$.

As an immediate consequence, it follows that the Group 8 from Table
\ref{t332} also has a finite-index subgroup that admits an epimorphism
onto a nonabelian free group.

\begin{verbatim}
gap> F:=FreeGroup(2);;
gap> x:=F.1;; y:=F.2;;
gap> W:=(x*y)^4 * (x*y^2)^3 * (x*y) * (x*y^2) * (x*y)^2 * (x*y^2);;
gap> G:=F/[x^2,y^3,W^2];;
gap> H:=Subgroup(G,[(G.1*G.2)^5]);;
gap> P:=PresentationNormalClosure(G,H);;
gap> gg:=GeneratorsOfPresentation(P);;
gap> for i in [1,3,7] do AddRelator(P,gg[i]); od;
gap> SimplifyPresentation(P);;
#I  there are 3 generators and 0 relators of total length 0
\end{verbatim}

3)
This GAP session shows that Group 7 in Table \ref{t332} is large,
following the proof in \cite{Will}.
$$G_7=\< x,y|x^3=y^3=(xyxyx^2y^2x^2yxy^2)^2=1\>.$$
It has a subgroup of index $12$ which admits an epimorphism onto
$\Z\ast\Z_2$.

\bigskip
\begin{verbatim}
gap> F:=FreeGroup(2);;
gap> x:=F.1;; y:=F.2;;
gap> W:=x*y*x*y*x^2*y^2*x^2*y*x*y^2;;
gap> G:=F/[x^3,y^3,W^2];;
gap> a:=G.1;; b:=G.2;;
gap> s1:=b*a^2;; s2:=a^2*b*a*b^2*a^2;;
gap> s3:=a*b*a*b^2*a^2*b*a;;
gap> s4:=a*b*a*b*a*b*a*b^2*a^2*b^2*a^2;;
gap> H:=Subgroup(G,[s1,s2,s3,s4]);;
gap> Index(G,H);
12
gap> P:=PresentationSubgroup(G,H);;
gap> gg:=GeneratorsOfPresentation(P);
[ _x1, _x2, _x3, _x4, _x5, _x6 ]
gap> AddRelator(P,gg[1]*gg[2]^-1*gg[1]*gg[2]^-1);
gap> AddRelator(P,gg[2]*gg[3]^-1);
gap> SimplifyPresentation(P);
#I  there are 2 generators and 1 relator of total length 4
gap> TzPrint(P);
#I  generators: [ _x1, _x2 ]
#I  relators:
#I  1.  4  [ -1, 2, -1, 2 ]
\end{verbatim}

4)
This GAP session uses the functions A5Poly,
A4A5Poly, S4A5Poly, and A4S5A5Poly (see the ancillary file {\em gtg232.g}
for function listings)
to compute all words in $\<x|x^2=1\>*\<y|y^3=1\>$, up to equivalence,
with trace polynomial of the form $(\l^2-1)^a(\l^2-2)^b(\l^4-3\l^2+1)^c$,
where $a,b\le 1$ and $c\le 3(a+b+1)$.
Words are output as even-length lists of positive integers:
$[a(1),b(1),\dots,a(t),b(t)]$ is shorthand for
$(xy)^{a(1)}(xy^2)^{b(1)}\cdots (xy)^{a(t)}(xy^2)^{b(t)}$.

The time requirement of each of these GAP functions grows
at least exponentially with the input parameter $c$.
For small values of $c$, the runtime is essentially
instantaneous.  But the final run (corresponding to the case $a=b=1,c=9$)
took close to 6 hours of CPU time on a 3GHz processor.  Thus it appears
that the theoretical limits supplied by Theorems \ref{tpoly} and \ref{tpoly2}
are not very far short of the practical limits for this impementation.

\begin{verbatim}
gap> A5Poly(1);
[ [ 3, 1 ] ]
gap> A5Poly(2);
[ [ 4, 1, 1, 2 ] ]
gap> A5Poly(3);
[ [ 4, 3, 1, 1, 2, 1 ] ]
gap> A4A5Poly(1);
[ [ 2, 1, 1, 2 ] ]
gap> A4A5Poly(2);
[ [ 3, 1, 1, 2, 1, 2 ] ]
gap> A4A5Poly(3);
[ [ 4, 2, 1, 1, 1, 2, 1, 2 ] ]
gap> A4A5Poly(4);
[ [ 4, 2, 1, 1, 2, 1, 1, 3, 1, 2 ] ]
gap> A4A5Poly(5);
[ [ 4, 2, 1, 2, 1, 3, 3, 1, 1, 2, 1, 1 ] ]
gap> A4A5Poly(6);
[  ]
gap> S4A5Poly(1);
[ [ 4, 2 ] ]
gap> S4A5Poly(2);
[ [ 4, 3, 2, 1 ] ]
gap> S4A5Poly(3);
[ [ 5, 3, 2, 1, 1, 2 ] ]
gap> S4A5Poly(4);
[ [ 4, 4, 2, 1, 1, 2, 3, 1 ] ]
gap> S4A5Poly(5);
[ [ 4, 4, 1, 1, 2, 3, 3, 1, 2, 1 ] ]
gap> S4A5Poly(6);
[  ]
gap> A4S4A5Poly(1);
[ [ 3, 2, 1, 2 ] ]
gap> A4S4A5Poly(2);
[  ]
gap> A4S4A5Poly(3);
[ [ 4, 2, 1, 1, 2, 3, 1, 2 ], [ 4, 3, 1, 2, 1, 1, 2, 2 ],
  [ 4, 3, 2, 2, 1, 2, 1, 1 ] ]
gap> A4S4A5Poly(4);
[ [ 4, 2, 1, 2, 3, 3, 1, 2, 1, 1 ], [ 4, 3, 1, 2, 1, 2, 3, 2, 1, 1 ] ]
gap> A4S4A5Poly(5);
[  ]
gap> A4S4A5Poly(6);
[ [ 4, 2, 1, 2, 1, 1, 3, 2, 1, 4, 3, 2, 1, 1 ],
  [ 4, 1, 1, 2, 3, 1, 1, 2, 1, 2, 3, 4, 1, 2 ],
  [ 4, 3, 1, 1, 2, 1, 2, 2, 3, 4, 1, 2, 1, 1 ],
  [ 4, 3, 1, 1, 2, 1, 3, 4, 1, 2, 1, 1, 2, 2 ] ]
gap> A4S4A5Poly(7);
[  ]
gap> A4S4A5Poly(8);
[  ]
gap> A4S4A5Poly(9);
[  ]
\end{verbatim}

\begin{table}[ht]
\caption{Words in $\Z_3\ast\Z_3$ with trace polynomial as in Theorem \ref{tpoly}}\label{t332}

\begin{center}
\begin{tabular}{|l||l|r|c|}
\hline
& $W(x,y)$ & $\tau(\l)$ & SCC\\ 
\hline
1 & $xy$ & $\l$ & NO\\
2 & $xyxy^2$ & $\l^2-\l-1$ & NO\\
3 & $xyx^2y^2$ & $\l(\l-1)$ & NO\\
4 & $xyxyx^2y^2$ & $\l(\l^2-\l-1)$ & NO\\
5 & $xyxyx^2yx^2y^2$ & $(\l^2-\l-1)^2$ & NO\\
6 & $xyxy^2x^2yx^2y^2$ & $\l(\l-1)(\l^2-\l-1)$ & NO\\
7 & $xyxyx^2y^2x^2yxy^2$ & $\l(\l^2-\l-1)^2$ & NO\\
8 & $xyxyx^2y^2x^2yx^2yxy^2$ & $(\l^2-\l-1)^3$ & NO\\
9 & $(xyxyx)(y^2x^2y^2x)(yx^2yx^2y^2)$ & $\l(\l^2-\l-1)^3$ & YES\\
10 & $(xyxy)(x^2y^2x^2yx)(y^2x^2yx^2y^2xy^2)$ & $\l(\l-1)(\l^2-\l-1)^3$ & YES\\
11 & $(xyxy)(x^2y^2x^2yx^2)(y^2xy^2xyx^2y^2)$ & $\l(\l-1)(\l^2-\l-1)^3$ & YES\\
12 & $(xyxy)(x^2y^2xy^2x^2y^2)(xyx^2yx^2y^2)$ & $\l(\l-1)(\l^2-\l-1)^3$ & YES\\
13 & $(xyxy)(x^2y^2x^2y^2)(xy^2x^2y^2xyx^2yx^2y^2)$ & $\l(\l^2-\l-1)^4$ & YES\\
14 & $(xyxy)(x^2y^2xy^2x^2yxy)(x^2y^2x^2yx^2y^2xy^2)$ & $\l(\l-1)(\l^2-\l-1)^4$ & YES\\
15 & $(xyxy)(x^2y^2x^2y^2)(xy^2x^2yx^2y^2x^2yxyx^2y^2xy^2)$ & $\l(\l^2-\l-1)^5$ & YES\\
16 & $(xyxyx^2y^2)(x^2yxy^2xy^2x^2y^2x^2)(yxy^2xyx^2yx^2y^2x^2yxy^2)$ & $\l(\l-1)(\l^2-\l-1)^6$ & YES\\
17 & $(xyxyx^2y^2x^2)(yxy^2xyx^2yx^2y^2x^2yxy^2x)(y^2x^2y^2x^2yxy^2)$ & $\l(\l-1)(\l^2-\l-1)^6$ & YES\\
18 & $(xyxy^2x^2yx)(yx^2y^2xy^2xyxy^2)(x^2y^2x^2yx^2yxy^2xyxy)$ & $\l(\l-1)(\l^2-\l-1)^6$ & YES\\
19 & $(xyx^2y^2x^2yx^2)(y^2xy^2xyxy^2x^2)(y^2x^2yxy^2x^2yx^2yxy^2xy)$ & $\l(\l-1)(\l^2-\l-1)^6$ & YES\\
\hline
\end{tabular}
~\\
\end{center}
The final column indicates whether or not $W$ satisfies the small-cancellation hypotheses of Corollary
\ref{sc33}.  In those cases where it does, the bracketing indicates a subdivision of $W$ into three
non-pieces of length $\ge 4$: $W\equiv U_1\cdot U_2\cdot U_3$.
\end{table}

\newpage
\begin{table}[ht]
\caption{Words in $\Z_2\ast\Z_3$ with trace polynomial as in Theorem \ref{tpoly2}}\label{t232a}
\centerline{(Numeration corresponds to related words in Table \ref{t332}.)}

\begin{center}
Part 1: Short words.  These groups are already known or can be easily analysed.
\begin{tabular}{|l||l|r|r|}
\hline
& $W(u,v)$ & $\tau(\l)$ & Size of $H$\\ 
\hline
1a & $uvuv$ & $\l^2-2$ & 24\\
1b & $uvuv^2$ & $\l^2-1$ & 24\\
2 & $uvuvuvuv^2$ & $\l^4-3\l^2+1$ & 120\\
3 & $uvuvuv^2uv^2$ & $(\l^2-1)(\l^2-2)$ & 576\\
4a & $uvuvuvuvuv^2uv^2$ & $(\l^2-2)(\l^4-3\l^2+1)$ & 2880\\
4b & $uvuvuv^2uv^2uvuv^2$ & $(\l^2-1)(\l^4-3\l^2+1)$ & large\\
5 & $uvuvuvuvuv^2uv^2uvuv^2$ & $(\l^4-3\l^2+1)^2$ & abelian-by-finite\\
6 & $uvuvuvuv^2uv^2uvuv^2uv^2$ & $(\l^2-1)(\l^2-2)(\l^4-3\l^2+1)$ & 424673280\\
8 & $uvuvuvuvuv^2uv^2uv^2uvuv^2uvuvuv^2$ & $(\l^4-3\l^2+1)^3$ & large\\
\hline
\end{tabular}
~\\
\end{center}

\begin{center}
Part 2: Small cancellation words.  Bracketing gives $W\equiv U_1\cdot U_2\cdot U_3$
as in Corollary \ref{sc23}.
\begin{tabular}{|l||l|}
\hline
& $W(u,v)$ \\
\hline
9b & $(uvuvuvuv)(uv^2uv^2uvuv^2uv^2)(uvuv^2uvuv^2uv^2)$ \\
10 & $(uvuvuvuv)(uv^2uv^2uv^2uvuv)(uv^2uv^2uvuv^2uv^2uvuv^2)$ \\
11 & $(uvuvuvuv)(uv^2uv^2uv^2uvuv^2)(uv^2uvuv^2uvuvuv^2uv^2)$ \\
13b & $(uvuvuvuv)(uv^2uv^2uvuv^2uv^2)(uv^2uvuv^2uvuvuv^2uvuv^2uv^2)$ \\
14a & $(uvuvuvuv)(uv^2uv^2uvuv^2uv^2uvuv)(uvuv^2uv^2uv^2uvuv^2uv^2uvuv^2)$ \\
14b & $(uvuvuvuv)(uv^2uv^2uv^2uvuv^2uv^2)(uvuv^2uv^2uvuvuvuv^2uv^2uvuv^2)$ \\
15b & $(uvuvuvuv)(uv^2uv^2uvuv^2uv^2uv)(uv^2uv^2uv^2uvuvuvuv^2uvuv^2uv^2uvuv^2)$ \\
16 & $(uvuvuvuvuv^2uv^2)(uv^2uvuvuv^2uvuv^2uv^2uv^2uv^2)(uvuvuv^2uvuvuv^2uvuv^2uv^2uv^2uvuvuv^2)$ \\
17 & $(uvuvuvuvuv^2uv^2uv^2)(uvuvuv^2uvuvuv^2uvuv^2uv^2uv^2uvuvuv^2uv)(uv^2uv^2uv^2uv^2uvuvuv^2)$ \\
18 & $(uv^2uv^2uv^2uvuvuv^2uv^2)(uvuv^2uv^2uvuv^2uvuvuvuv^2)(uv^2uv^2uv^2uvuv^2uvuvuv^2uvuvuvuv)$ \\
19 & $(uvuv^2uvuvuv^2uvuvuvuv)(uv^2uv^2uv^2uvuv^2uv^2uvuv^2uvuvuv)(uv^2uv^2uv^2uv^2uvuvuv^2uv^2)$ \\
\hline
\end{tabular}
~\\
\end{center}

\begin{center}
Part 3.  Cases remaining open
\begin{tabular}{|l||l|}
\hline
& $W(u,v)$\\
\hline
7a & $uvuvuvuvuv^2uv^2uv^2uvuvuv^2$ \\
7b & $uvuvuvuv^2uv^2uvuv^2uv^2uvuv^2$ \\
9a & $uvuvuvuvuvuv^2uv^2uv^2uvuvuv^2uvuv^2uv^2$ \\
12 & $uvuvuvuvuv^2vu^2uvuv^2uv^2uv^2uvuvuv^2uvuv^2uv^2$ \\
13a & $uvuvuvuvuv^2uv^2uv^2uv^2uvu^2v^2uv^2uvuvuv^2uvuv^2uv^2$ \\
15a & $uvuvuvuvuv^2uv^2uv^2uv^2uvuv^2uv^2uvuv^2uv^2uv^2uvuvuvuv^2uv^2uvuv^2$ \\
\hline
\end{tabular}
~\\
\end{center}

\end{table}

\end{document}